\documentclass[11pt, reqno ]{amsart}
\usepackage{geometry}                
\usepackage{graphicx}
\usepackage{amssymb}
\usepackage{epstopdf}
\usepackage{amsmath}
\usepackage{mathrsfs}
\usepackage[all]{xy}
\usepackage{color}
\usepackage{tikz}
\usepackage[colorlinks=true, pdfstartview=FitV,linkcolor=blue,citecolor=blue,urlcolor=blue]{hyperref}

    \advance \topmargin by -\headheight
    \advance \topmargin by -\headsep

    \evensidemargin \oddsidemargin
\numberwithin{equation}{section}

\begin{document}

 \newtheorem{theorem}{Theorem}[section]
  \newtheorem{prop}[theorem]{Proposition}
  \newtheorem{cor}[theorem]{Corollary}
  \newtheorem{lemma}[theorem]{Lemma}
  \newtheorem{defn}[theorem]{Definition}
  \newtheorem{ex}[theorem]{Example}
    \newtheorem{conj}[theorem]{Conjecture}

 \newcommand{\cx}{{\bf C}}
\newcommand{\la}{\langle}
\newcommand{\ra}{\rangle}
\newcommand{\res}{{\rm Res}}
\newcommand{\expp}{{\rm exp}}
\newcommand{\Sp}{{\rm Sp}}
\newcommand{\lgp}{\widehat{\rm G}( { F} ((t)) )}
\newcommand{\svee}{\scriptsize \vee}
\newcommand{\deff}{\stackrel{\rm def}{=}}
\newcommand{\cal}{\mathcal}

\title{A Construction of Representations of Loop Group and Affine Lie Algebra of $\frak {sl}_n$     }

\author[Xuanzhong Dai, Corresponding Author: Yongchang Zhu]{Xuanzhong Dai,   Yongchang Zhu}
\address{Department of Mathematics, The Hong Kong University of Science and Technology, Clear Water Bay, Kowloon, Hong Kong}
\email{mazhu@ust.hk, xdaiac@connect.ust.hk    }

\thanks{The research is supported by Hong Kong RGC grant  16305715.}

\maketitle

\maketitle

\section { Introduction }

The loop groups of simple Lie groups and their algebraic counterparts affine Kac-Moody algebras
 are infinite dimensional generalizations of simple Lie groups and simple Lie algebras.
  Many results of latter objects have infinite dimensional generalizations. In this work we construct
  a family of representations of certain central extension of loop group $ SL_n ( {\Bbb R} (( t))) $
   and its affine Lie algebra $\widehat {\frak {sl}}_n $.  This family is parametrized by the characters of the  multiplicative group
  \begin{equation}
      G \deff \{   a( t) \in {\Bbb R} (( t))^* \; | \;  {\rm the \; leading \; coefficient \; of}
       \; a ( t) \; { \rm is \; } \pm 1 \}.
  \end{equation}
    The group  $G$ is a subgroup of a certain central extension of loop $GL_1 ( {\Bbb R} (( t))) $.
      Our construction can be interpreted loosely as  an affine analog of  the theta correspondence for the dual pair $GL_1 $ and $GL_n$.
      In \cite{Z}, the second author constructed
  a Weil representation for metaplectic loop group of $SP_{2n}$, in which  the above group $G $ and the lifting of
   $SL_n ({\Bbb R} (( t))) $ form a commuting pair.  Our construction can be viewed as the
    theta lifting of characters of $G$ to representations of
    the central extension of   $SL_n ({\Bbb R} (( t))) $.
 Since this  connection is formal and the calculations are rather involved, we choose
   to present our construction in a way independent of the work in \cite{Z}.

   The underline space
of our  representation consists of  functions  on $ { \Bbb R }^n (( t ))  $ that is periodic under the translation by elements in
 $ {\Bbb R}^n [[ t ]] \subset  { \Bbb R }^n (( t )) $, so it is a certain subspace of functions on $    { \Bbb R }^n (( t )) /   {\Bbb R}^n [[ t ]] = {\Bbb R}^n [ t^{-1} ] t^{-1} $.
   The operators appear in our representation are integral operators.
  The operators for the corresponding representation of affine Lie algebras are infinite sums of quadratic differential operators,
   which is similar to early works of free field realizations \cite{Fr} \cite{KP} \cite{FF}, but the operator for a negative mode
     element in $\widehat{\frak {sl}}_n $ has  one term as an integral operator. And the representations we obtain
      are not highest weight modules or in the category ${\cal O}$.

 To introduce our construction, we begin with two similar but much simpler cases.   Let $ C ( {\Bbb R} / {\Bbb Z}) $ be the space of continuous functions on the circle $ {\Bbb R} / {\Bbb Z} $.
 Although the multiplicative semigroup  ${\Bbb Z}_{\ne 0 } $ doesn't act on ${\Bbb R} / {\Bbb Z}$ in any natural way, it acts on the function space  $ C ( {\Bbb R} / {\Bbb Z}) $ as Hecke operators:
 for a non-zero integer $ n $, let $ \pi(n) $ be the operator given by
 \begin{equation}\label{1.1}
  \pi ( n )  f ( x ) = \sum_{ i =1}^n f ( \frac x n + \frac i  n ).
\end{equation}
  It is easy to check that
  $ \pi ( m ) \pi ( n ) = \pi ( mn ) $.
  The common eigenfunctions for $\pi ( n )$'s with certain continuity properties are classified in \cite{M}.

  The second case is related to representations of $GL_n ( F) $, where $F$ is a non-Archmedean local field. Let $R $ be the ring of integers of $F$. We consider the space
  $ Fun (  F^n / R^n ) $ of complex valued functions on $F^n $ that is periodic with period $R^n $. Let $ M_n' ( R)  $ be the semi-group of $n \times n $-matrices with entries in $R$ and with non-zero determinant.
  For $ g \in  M_n' ( R)  $ and $ f(x) \in Fun (  F^n / R^n ) $, we view $f$ as a function on $F^n$ which is invariant under the translations by $R^n$. The function
   $ f ( g^{-1} x ) $ is periodic with periods $ g  R^n $.  Since $ g $ has entries in $R$, $ g R^n \subset R^n$, and the quotient space $  R^n /    g R^n $ is a finite set.
   To get a periodic function with period $ R^n$, we take the average of  $ f ( g^{-1} x ) $ over    $  R^n /    g R^n $, we get linear  operator $ \pi( g ) :   Fun (  F^n / R^n ) \to  Fun (  F^n / R^n )     $ :
     \begin{equation}\label{1.2}
  \pi ( g )  f ( x ) = \sum_{ r \in R^n / g R^n } f (  g^{-1}  ( x  + r )  ) =   \sum_{ r \in  g^{-1} R^n /  R^n } f (  g^{-1}   x  +  r  )         .
   \end{equation}
  It is easy to prove that $ \pi ( g_1 ) \pi ( g_2 ) = \pi ( g_1 g_2 )$. So we have a representation of semi-group $ M_n' ( R)  $ on  $ Fun (  F^n / R^n ) $.
   To get a representation of $ GL_n ( F)$, we invert the operator $ \pi ( \beta I_n )$, where $\beta \in R $ is a prime. A direct way is to consider the space
   \[  V_{\lambda } =  \{   f \in   Fun (  F^n / R^n ) \; | \;   \pi ( \beta I_n ) f = \lambda f \}. \]
 For  $ g \in GL_n ( F)$,  we decompose it as $  g = \beta^k g' $, where $k \in {\Bbb Z}$, $ g' \in M_n' ( R ) $, then
  $ \pi_\lambda  ( g ) = \lambda^k \pi ( g' )$ is independent of the decomposition, and we have $ \pi_\lambda ( g_1 g _2 ) = \pi_\lambda ( g_1 ) \pi_\lambda ( g_2 ) $.
  In general, for a rational morphism of a reductive groups $ G \to GL_n $, then $V_\lambda $ is a representation of $G$ via the pull-back. In particular, we have an embedding
   $ GL_n \times GL_n \to GL_{n^2} $ by the action $ ( g_1 , g_2 ) x = g_1 x g_2^t $, where $x$ is an $n\times n$ matrix.  Our representation
    $V_\lambda $ of   $ GL_n \times GL_n $ can be used to give an interpretation of gamma factors of irreducible admissible representations of $GL_n ( F)$ as defined in \cite{GJ}. We will discuss this connection in a separate work.

\

  The field $ {\Bbb R} (( t)) $ and a non-Archedemean local field have many similarities.
  The above construction generalizes to $SL_n $ over  $ {\Bbb R} (( t)) $. From now on,  we  denote $F$ the field $ {\Bbb R} ((t))$, and $R$ the subring $ {\Bbb R} [[ t]] $. Let
  $M_n' ( R ) = M_n' ( {\Bbb R} [[ t]]  ) $ be the semigroup of $n\times n$ matrices over $R$ with non-zero determinant.
  We follow the construction (\ref{1.2}). For a function $f$ on  $ Fun (  F^n / R^n )$, i.e., $f$ is a function on $F^n$ with periods $R^n$,  and $ g \in   M_n' ( R ) $.
  $ f ( g^{-1} x )$ is periodic with periods $ g R^n $. But now $ R^n / gR^n \cong g^{-1} R^n / R^n   $ is not finite, it is a finite dimensional vector space over ${\Bbb R}$. To generalize (\ref{1.2}), we
   choose a Haar measure $\mu $ on $ g^{-1} R^n / R^n   $, we introduce operator $ \pi ( g , \mu ) $ by
   \begin{equation}\label{1.3}
  \pi ( g , \mu )f(x) \deff    \int_{   g^{-1} R^n / R^n }  f (  g^{-1}   x  +  v  ) \mu ( d v ).
  \end{equation}
  The  collection of all pairs $( g , \mu )$ form a semigroup $ \widehat {M_n} ( R )' $
  that is a central extension of $ M_n' ( R )$.  We can again consider the eigenspace of $ \pi ( t I_n , \mu )$. Now the situation is different from the local field case in that
   $ \pi ( t I_n , \mu )$ don't commutes with $ \pi ( g , \mu )$ for general $ g \in M_n' ( R) $, it only commutes with operators
    $ \pi ( g , \mu )$ with $ {\rm det } \, g $ has leading coefficient $\pm 1 $, as a result, the eigenspace is a representation of a subgroup of
     $ \widehat { GL_n }( F)  $ that contains   $ \widehat { SL_n }( F)  $.
     We will study various properties of this representation,
   and derive a formula for  the corresponding  representation of affine Kac-Moody algebra of $\frak {sl}_n$.

 \

  The structure of this paper is as follows.  In Section 2, we will give details of the construction of the representation of loop group $\widehat{SL_n} ( {\Bbb R} (( t)))$
  mentioned above and study its properties.
  In Section 3, we define a ``maximal compact subgroup" for $\widehat{SL_n}(\Bbb R)$ and construct a Gaussian function fixed by it.  In Section 4, we derive the formula of the corresponding representation of the affine Kac-Moody algebra $\widehat {\frak {sl}_n}$, and this representation has level $1$.
  In Section 5, we consider the dual action of $\widehat{\frak {sl}_n}$ and
  construct various highest weight representations of level $-1$.
   In Section 6. we construct a representation of $\widehat  {\frak {sl}_n} \oplus  \widehat  {\frak {sl}_n}$
  via embedding $ GL_n \times GL_n \to GL_{n^2}$, the representation has the special property that the level is critical.
    And we construct Whittaker functionals in the representation.

   We wish to thank I.Frenkel for discussions.

\

\section { A Construction of  a Representation of Loop $SL_n$.}\label{s2}

We continue to use the notations
 \begin{equation}\label{2.1}
  F = {\Bbb R}(( t)) , \; \; \;  \; \; R = {\Bbb R} [[ t ]], \; \;  \;  \; \; R_- = {\Bbb R} [ t^{-1}] t^{-1}, \; \; \; \; F = R \oplus R_-.
 \end{equation}

Recall that a function $f$ on a finite dimensional vector space over ${\Bbb R}$ is called a Schwartz function
 if  all its partial derivatives $ \partial^I f $  of  arbitrary order is rapidly decay in the sense that for all $m$,
  $  ( 1 + | x |^m ) \partial^I f ( x ) $ is bounded.
  A key property of the space  ${\cal S} ( V)$ of all Schwartz functions on $V$ is that it is closed under
  taking derivatives, multiplying by polynomials and taking partial Fourier transforms.

For an infinite dimensional vector space $V$ over ${\Bbb R}$ such as $ V = R_-^n  = {\Bbb R}^n [ t^{-1} ] t^{-1} $,
a function $ f $ on $V$ is called a
Schwartz function if the restriction of $f$ to every finite dimensional space is a Schwartz function.
We denote by  ${\cal S} ( F^n / R^n  ) = {\cal S} ( R_-^n )$ the space of all Schwartz functions on $R_-^n = {\Bbb R}^n [ t^{-1} ] t^{-1} $.

Let $M_n' ( R )    $ be the semi-group of $n\times n$ matrices over $R$ with non-zero determinant. We first construct
   a representation of a central extension of $M_n' ( R )    $  on  ${\cal S} ( F^n / R^n  )$.

Let $ \pi_- : F^n \to  R_-^n $ be the projection map with respect to the decomposition $ F^n = R^n \oplus  R_-^n$.
A function  $ f \in S( R_-^n ) $ gives a function on $ F^n$, which we still denote by $f$, by the pull-back along $\pi_- $, so
$ f ( x )= f ( \pi_- ( x ) )$.  It is clear that $ f $ is $ R^n$-periodic. Conversely, $ R^n$-periodic on $F^n$ is the pull-back
 for a unique function on $ R_-^n$. We will not distinguish the $R^n$-periodic functions on $F^n$ and functions on $R_-^n$.

For each $ g \in  M_n' ( R )$ and an $R^n$-periodic function $f$ on $F^n$,  $ f ( g^{-1} x ) $ is $ g R^n$-periodic.
  Since $g$ has entries in $ {\Bbb R} [[ t]] $, we have $   g  R^n  \subset R^n$, therefore
  $R^n  \subset g^{-1} R^n $.
   The quotient space  $ g^{-1} R^n   / R^n  $ is a finite dimensional
   space over ${\Bbb R}$.
  If $  {\rm det}\, g = t^N ( c_0 + c_1 t + \cdots ) $, $c_0 \ne 0 $, then
  \[ {\rm dim}_{\Bbb R} \,   g^{-1}  R^n  / R^n =   N .\]
  We denote
   \[   V_g =   g^{-1}    R^n  / R^n  .\]
 Let $\mu$ be a Haar measure on $V_g$,  then we consider the operator
  $\pi ( g , \mu )$ given by
  \begin{equation}\label{2.2}  \pi(g,\mu) f ( x ) : = \int_{ V_g } f ( g^{-1}  x + y ) \mu(dy)  .\end{equation}
 For $ a \in R^n $,
 \[ \pi(g,\mu) f ( x +a  )= \int_{ V_g } f ( g^{-1}  x + g^{-1} a + y ) \mu(dy)\]
 since $ g^{-1}a \in g^{-1} R^n $, we change variable  $ g^{-1} a + y \to y$, using the translation variance of the measure $ \mu$,
  we see that
\[ \int_{ V_g } f ( g^{-1}  x + g^{-1} a + y ) \mu(dy) = \int_{ V_g } f ( g^{-1}  x + y ) \mu(dy) = \pi(g,\mu) f ( x   )
\]
 So $\pi(g,\mu) f ( x   )$ is $R^n$-periodic. Since a linear change of variable of a Schwartz function is again  a Schwartz function
  so $  f ( g^{-1}  x )$ is a Schwartz function.   And  if $ h ( x , y )$ is a Schwartz function on $ ( x , y) \in U \times V$
   for two finite dimensional spaces $U $ and $V$, then $ \int_V h ( x , y ) dy $ is a Schwartz function on $  x\in U$ for any Haar
    measure $dy$ on $V$, we see that $\pi(g,\mu) f ( x ) \in {\cal S} ( R_-^n)$.

\

Let $\widehat{M_n} ( R )'$ be the set of pairs $(g,\mu)$, where $g $ is an element in $M_n  ( R )'$, and $\mu$ is a Haar measure on $V_g$.
 We define a multiplication on $\widehat{M_n}(R)'$ so that $ \pi ( g , \mu )$ is a representation on ${\cal S} ( R_-^n)$.
For $(g_1,\mu_1), (g_2,\mu_2) \in \widehat{M_n}(R)'  $,
 \begin{equation}\label{2.3} (g_1,\mu_1)  (g_2,\mu_2)= (g_1g_2,\mu_1\ast \mu_2).\end{equation}
 where $\mu_1\ast \mu_2$ is the Haar measure on $ V_{ g_1 g _2 } $ given as follows.
  The chain  $ g_1 g_2  R^n  \subset g_1  R^n  \subset   R^n $
   induces the chain $   R^n  \subset g_2^{-1}  R^n  \subset g_2^{-1}  g_1^{-1} R^n $.
  Let  $\mu_1'$ be the  pushforward measure on $\tilde V = g_2^{-1}g_1^{-1} R^n  / g^{-1}_2 R^n$ of  $\mu_1$ under the
   isomorphism $\phi: V_{g_1} = g_1^{-1} R^n  / R^n \to g_2^{-1}g_1^{-1} R^n  / g^{-1}_2 R^n $,
   sending $y\;  ({\rm mod} \, R^n)  $ to $g_2^{-1}y \; ({\rm mod} \, g^{-1}_2 R^n) $.
  Consider the
short exact sequence:
\[
0\to g_{2}^{-1} R^n  / R^n  \to g_2^{-1} g_1^{-1} R^n / R^n \to g^{-1}_2g^{-1}_1 R^n  / g_2^{-1} R^n \to 0
\]
where the second arrow refers to the inclusion map, and the third arrow refers to the corresponding quotient map.
   Then we define the Haar measure $\mu_1\ast \mu_2$ on  $V_{g_1g_2} =  g_2^{-1} g_1^{-1} R^n / R^n$
  such that,  for $f\in L^1 ( g_2^{-1} g_1^{-1} R^n / R^n   )$,
\begin{eqnarray}\label{2.4}
\int_{V_{g_1g_2}} f(x) (\mu_1\ast \mu_2)(dx) &=& \int _{\tilde V} \int_{V_{g_2}} f(\tilde y+z)\mu_2(dz)\mu_1'(d\tilde y) \nonumber \\
&=& \int _{ V_{g_1}} \int_{V_{g_2}} f( g_2^{-1} y+z)\mu_2(dz)\mu_1(d y)
\end{eqnarray}

\begin{prop} \label{proposition 2.1} $\widehat{M_n}( R )'$ is a semi-group under the product (\ref{2.3}).
(\ref{2.2}) gives a representation of  $\widehat{M_n}( R )'$ on ${\cal S} ( R_-^n )$.
\end{prop}

\noindent {\it Proof.}  It is clear that $ V_{ I_n} = \{ 0 \}$, and let $ \mu_{st} $ denote the counting measure on  $ V_{ I_n}$.
Obviously $(I_n,\mu_{st})\in \widehat{M_n}( R )'$   is the identity element .
To prove (\ref{2.3}) is associative for $(g_1,\mu_1),(g_2,\mu_2),(g_3,\mu_3)$,
let $f$ be arbitrary  integrable function  on $  V_{ g_1 g_2 g_3 }  = g_3^{-1}g_2^{-1}g_1^{-1} R^n / R^n $,
 from the definition (\ref{2.4}) of $\ast$, we have
\begin{align*}
&\int_{V_{g_1g_2g_3}} f(v)( \mu_1\ast (\mu_2\ast \mu_3)) (dv)\\
 = &\int_{V_{g_1}} \int_{V_{g_2g_3}} f((g_2g_3)^{-1} y+z ) (\mu_2\ast\mu_3)(dz)\mu_1(dy)\\
=& \int_{V_{g_1}}\int_{V_{g_{2}}}\int_{V_{g_3}} f((g_2g_3)^{-1}y +g_3^{-1}z + w  ) \mu_{3}(dw) \mu_{2}(dz) \mu_{1}(dy)
\end{align*}
Similarly  we can show $\int_{V_{g_1g_2g_3}} f(v)( (\mu_1\ast \mu_2)\ast \mu_3) (dv)$ also equals
 to the right hand side the above formula. Therefore the associativity is proved. Next we prove
$ \pi(g_1,\mu_1) \pi(g_2,\mu_2)=\pi(g_1g_2,\mu_1\ast \mu_2) $.
 For any $f\in {\cal S } (R_-^n )$,  we have
\begin{align*}
\pi(g_1g_2,\mu_1\ast \mu_2) f(x)=&  \int _{V_{g_1g_2}} f((g_1g_2)^{-1}x+y ) (\mu_1\ast \mu_2) (d y)\\
                                                    =& \int _{V_{g_1} }\int_{V_{g_2}} f((g_1g_2)^{-1}x+g_2^{-1}y+z) \mu_2(dz)\mu_1(dy)
                                                 \end{align*}
\begin{align*}
\pi(g_1,\mu_1)(\pi(g_2,\mu_2) f )(x)=&\int _{V_{g_1}} \pi(g_2,\mu_2) f ( g_1^{-1} x  + y    ) \mu_1 ( y ) \\
=&\int _{V_{g_1}} \int_{ V_{g_2}}  f ( g_2^{-1} ( g_1^{-1}  x  + y ) + z    ) \mu_2 ( d z ) \mu_1 ( y )
 \end{align*}
This completes the proof.
\qed

\

  Since every element $ g \in GL_n (F ) $ can be written as $ g =  t^k h $ for $ h \in M_n ( R ) ' $ and $ k \leq 0 $, we need to
   invert the operator $ \pi ( t I_n , \mu )$ to extend the $\pi$-action.
   A convenient way is to consider the eigenspace of $  \pi ( t I_n , \mu )$, and we then are lead to study
  the elements in $\widehat{M_n}( R )'$ that commutes with $ \pi ( t I_n , \mu )$.

  \

  For $ u \in GL_n (  R  ) $, $ u R^n   = R^n $,  so $V_u = \{ 0 \}$.
   We denote $\mu_{st}$ for the counting measure on the one point set $V_u = \{ 0 \} $. Then for any $f\in \mathcal S( R_-^n ) $, we have
   \begin{equation}\label{2.5}  \pi(u,\mu_{st} )f ( x )  =  f ( u^{-1}  x )  .\end{equation}
 For $ {\bf k} = ( k_1 , \cdots , k_n ) \in {\Bbb Z}_{\geq 0 }^n $, let
    \[  t^{\bf k} = {\rm diag } ( t^{k_1} , \cdots , t^{k_m} ) \]
    It is easy to see that
 \[ V_{ t^{\bf k} } = {\rm Span} (    t^{-i_j } e_j  \, | \,  1\leq  j \leq m ,   1\leq  i_j  \leq k_j ), \]
 Let   $\mu_{ st }$ denote the Haar measure induced by the inner product on $V_{ t^{\bf k} }$ on which $\{   t^{-i_j } e_j \}$ is an
    orthonormal basis.

It is easy to show that,  for $ u, v \in GL_n ( R  )$,
\begin{equation}\label{2.6}  (uv,  \mu_{\rm st}  ) = (u,\mu_{\rm st})  (v,\mu_{\rm st}),\end{equation}
and
\begin{equation}\label{2.7}  (t^{\bf k + \bf k'},\mu_{\rm st}) =(t^{\bf k},\mu_{\rm st})  ( t^ {\bf k'},\mu_{\rm st})    .\end{equation}
We also have for $ c \in {\Bbb R} , c \ne 0 $,
\begin{equation}\label{2.8}
      ( t^k I_n , \mu_{st} )  ( c I_n , \mu_{st} )   = |c|^{k n }  ( c I_n , \mu_{st} ) ( t^k I_n , \mu_{st} ) .\end{equation}
More generally, for $ u \in GL_n ( {\Bbb R} [[ t]] ) $,
\begin{equation}\label{2.9}
     ( t^k I_n , \mu_{st} ) ( u , \mu_{st} )  = | {\rm det} (  u|_{t=0} ) |^k  ( u , \mu_{st} ) ( t^k I_n , \mu_{st} ).
      \end{equation}
where for $ a > 0 $, $ ( g , \mu ) \in  \widehat{M_n}( R )'$,  $ a ( g , \mu ) = ( g , a \mu )$ and  $ a \mu $ is the Haar measure
 on $V_g$ obtained by rescaling $\mu $ by $a$.
We prove (\ref{2.8}) by using the representation $ {\cal S} ( R_-^n )$, the other identities can be proved similarly.
\newline  {\it Proof}  of (2.8).
By (\ref{2.7}), it suffices to prove for $k=1$. For any $f\in \mathcal S ( R^n_- )$,
\begin{eqnarray*}
 \pi ( t I_n , \mu_{\rm st} )  \, \pi ( c I_n , \mu_{st} ) f ( x ) &=& \int_{ V_t}     \pi ( c I_n , \mu_{st} ) f (  t^{-1} x + y )\mu_{\rm st}( d y) \\
 &=&   \int_{ V_t}     f (  t^{-1} c^{-1} x + c^{-1} y ) \mu_{\rm st}(d y) \\
 &=&   |c|^n  \int_{ V_t}     f (  t^{-1} c^{-1} x +  y ) \mu_{\rm st}(d y) \\
 &=&   |c|^n \,  \pi ( c I_n , \mu_{st} )  \, \pi ( t I_n , \mu_{\rm st} )   f ( x ) .
\end{eqnarray*}

 We will write
   \[  \pi_{ t^k} =  \pi ( t^k I_n ,\mu_{\rm st}) , \; \; \; \; \;  \pi_c = \pi( c I_n,\mu_{\rm st} ).\]

  \

  For an eigenvalue $  \lambda \in {\Bbb C} $, $\lambda \ne 0 $,  we consider the eigenspace of $ \pi_t = \pi ( t I_n , \mu_{st} )  $
\begin{equation}
  {\cal S}(R_-^n )_\lambda =  {\cal S}(\mathbb R^n [t^{-1}]t^{-1})_\lambda = \{ f \in S( R_-^n )\; | \; \pi_{ t} f  = \lambda f  \}.
\end{equation}
We will prove this space is a representation of a central extension of a subgroup of $GL_n (F) $
that is slightly larger than $ SL_n ( F) $, which is given by
 \[  GL_n (F)_0 = \{ g \in  GL_n( F) \, | \; {\rm det} \, g \text{ has the leading coefficient } \pm 1   \} \]
where the leading coefficient refers to the coefficient of the lowest power of $t$.
Obviously $GL_n (F)_0$ is a subgroup of $GL_n (F)$.

\

 \begin{lemma}\label{lemma 2.2}
  For $ g \in  GL_n (F)_0\cap M_n ( R)' $, $ k \in {\Bbb Z}_{> 0 } $,
  and any Haar measure $\mu$ on $V_g$,  $ ( g , \mu ) $ and $ ( t^k I_n , \mu_{st} )$ commutes
   in  $\widehat{M_n}( R )'$, i.e.,
   \[  (g,\mu) ( t^k I_n , \mu_{st} )  =  ( t^k I_n , \mu_{st} ) (g,\mu).
  \]
\end{lemma}

\noindent {\it Proof.}
Since $ \pi_{t^k} = \pi_t^k $, it is enough to prove the case $k=1$.
By the Bruhat decomposition (see e.g., \cite{GR})
 \[  GL_n ( {\Bbb R} (( t))) = \cup   GL_n ( {\Bbb R} [[ t]])\, {\rm diag} ( t^{k_1} , \cdots , t^{k_n} ) \, GL_n ( {\Bbb R} [[ t]]), \]
 we can prove that
 \[  GL_n (F)_0 = \cup   GL_n (R)_0 \, {\rm diag} ( t^{k_1} , \cdots , t^{k_n} ) \, GL_n ( R)_0, \]
where
 \[  GL_n (R )_0 =  GL_n (F)_0 \cap   GL_n ( R ) =
 \{ g \in GL_n (R ) \, | \,  {\rm det} ( g |_{t=0} ) = \pm 1 \} .\]
We write $ g $ as $  g = h_1 t^{\bf k} h_2 $ according to the last decomposition,
$ {\bf k } = ( k_1 , \dots , k_n )$,  $ t^{\bf k} =  h_1^{-1} g h_2^{-1}$ has
 entries in $ R$, so all $ k_i \geq 0 $.  $  ( g , \mu )  =  ( h_1 , \mu_{st} )  ( t^{\bf k } , \mu_{st} )  ( h_2 , \mu_{st} ) $ up to a scalar,
the results follows from (\ref{2.7}) and (\ref{2.9}).
\hfill $\Box $

\

By Lemma \ref{lemma 2.2}, the space ${\cal S} (R_-^n )_\lambda$ is stable under $\pi(g,\mu)$ for
$g\in GL_n(F)_0 \cap M_n(R)'$.
  We may extend the $ \pi$-action on ${\cal S} (R_-^n )_\lambda$ to a central extension of the  group   $ GL_n (F)_0$.
  We first define the central extension.  Let  $\widehat{M_n}( R )_0'$ denote the set of pairs $ ( g , \mu ) \in \widehat{M_n}( R )'$ with
   $ g \in   GL_n ( F)_0 \cap M_n ( R )' $, it is a semisubgroup of  $ \widehat{M_n}( R )' $.
  Consider the direct product semigroups
   \[    \tilde G =   \{ t^k I \, | k \in {\Bbb Z} \} \times  \widehat{M_n}( R )_0', \]
  we define an equivalence relation $\sim $ on $\tilde G $ induced by
   \begin{equation}  ( t^{k +l } ,  \alpha ) \sim ( t^k ,    ( t^l I_n , \mu_{st} ) \alpha ),  \; \; l \in {\Bbb Z}_{\geq 0} . \end{equation}
  It is clear that $\sim $ is compatible with semi-group structure of $ \tilde G$, therefore the
  set $\tilde G / \sim $ of equivalence classes is a semigroup. We prove it is a group and is a central extension of $ GL_n (F )_0 $.

   \begin{lemma}\label{lemma 2.3}   $ \tilde G/ \sim $ is a group and    $ \tilde G/ \sim \to  GL_n ( F)_0 $, $ ( t^k I_n , ( g , \mu ) ) \mapsto
     t^k g $ is a surjective group homomorphism with kernel isomorphic to ${\Bbb R}_{>0} $.
   \end{lemma}

\noindent {\it Proof.}
The only non-trivial part is the the existence of inverse in  $ \tilde G/ \sim $. For $ \alpha \in G / \sim $ represented by  $ ( t^k I_n , ( g , \mu ) ) $,
 let $l $ large enough so $ t^l g^{-1} \in M_n ( R )'_0 $, choose an Haar measure $\nu $ on $ V_{ t^l g^{-1} }$, then
  \begin{equation}\label{2.12}     ( t^k I_n , ( g , \mu ) ) ( t^{ -k - l } I_n ,  (  t^l g^{-1} , \nu ) ) = ( t^{-l} , ( t^l I_n , \mu * \nu ) )  \end{equation}
$ ( t^l I_n , \mu * \nu ) = c  ( t^l I_n , \mu_{st}  ) $ for some $ c > 0 $, so the right hand side of (\ref{2.12}) is
\[    ( t^{ - l } I_n ,  (  t^l I_n  , \mu_{st}  ) c  ) \sim ( I_n , c ) \]
 which is invertible in $ \tilde G / \sim $.\qed

  We will the denote the group  $ \tilde G / \sim $ in the above lemma by $ \widehat {GL_n} (F)_0 $, which is a central extension
   of ${GL_n} (F)_0 $ by ${\Bbb R}_{>0} $:
   \begin{equation}
     1\to    {GL_n} (F)_0  \to   \widehat {GL_n} (F)_0 \to {\Bbb R}_{>0} \to 1 .
  \end{equation}
  And we denote
   the inverse image of $SL_n ({\Bbb R} (( t)) ) $ by  $ \widehat {SL_n} ( {\Bbb R} (( t)) ) $.

\begin{theorem}
  $ {\cal  S} (R^n_-)_\lambda$ is a representation of  $ \widehat {GL_n} ( F)_0 $,
   an element in $ \widehat {GL_n} ( F)_0 $ represented by $ ( t^l , ( g , \mu ) ) $ acts as
   \[    \pi ( t^l , ( g , \mu ) ) f ( x ) = \lambda^l \pi ( g , \mu ) f ( x )
    =  \lambda^l \int_{ V_g }  f ( g^{-1} x +  y ) \mu ( d y ) .\]
 And $ {\cal  S} (R^n_-)_\lambda$ is a representation of $ \widehat {SL_n} (F) $ by restriction.
\end{theorem}

\

\begin{lemma}\label{lemma 2.5}
On $\mathcal S(R^n_-)$, for any nonzero constant $c$, we have
 \begin{equation}\label{2.13}
  \pi_{c} \, \mathcal S(R^n_-)_\lambda  = \mathcal S(R^n_-)_{ |c|^n \lambda}
\end{equation}
 And $\pi_c$  is an isomorphism of the $ \widehat {SL_n} (F) $-modules
\end{lemma}

\noindent{\it Proof.}  From (\ref{2.8}) with $k=1$, we see that
 \[  \pi_{c} \, \mathcal S(R^n_-)_\lambda  \subset \mathcal S(R^n_-)_{ |c|^n \lambda} , \]
then using $ \pi_{c^{-1}} $, we see the equality of (\ref{2.13}) holds.
It remains to prove that
     $  \pi_{c} \pi(g) = \pi(g) \pi_c $ for $g \in \widehat {SL_{n}} ( F)$.
      Use the Bruhat decomposition of
       $\widehat {SL_{n}} ( F)$, it is enough to prove the cases $g \in SL_{n} (R) $ and
        $ g = ( t^{ l }   I_n ,  ( t^{\bf k }  , \mu_{st } ) ) $, where
    \[    t^{\bf k }   =  {\rm diag} ( t^{k_1} , \cdots , t^{k_n} ) \]  with $ k_1 + \cdots + k_n + l n = 0  $.
   The case $g\in SL_n (R)$ follows from (\ref{2.6}).
     For the later case,  let $f\in {\cal S} (R^n_-)_{\lambda } $,
     from (\ref{2.13}), $\pi_cf \in {\cal S} (R^n_-)_{|c|^n \lambda}$, we have,
      ${\rm dim}\,  V_{ t^{\bf k } } =  k_1 + \dots + k_n = - n l $,
       \begin{eqnarray*}
      \pi( t^{ l }   I_n ,  ( t^{\bf k }  , \mu_{st } ) ) \pi_c f ( x ) &=& (|c|^{ n  } \lambda )^l    \pi ( t^{\bf k }  , \mu_{st } )  \pi_c f  \\
       &=&  (|c|^{ n  } \lambda )^l   \int_{ V_{ t^{\bf k } } }  f ( c^{-1} t^{ - \bf k } x + c^{-1} y )   \mu_{st } ( d y ) \\
       &=&  \lambda ^l   \int_{ V_{ t^{\bf k } } }  f ( c^{-1} t^{ - \bf k } x +  y )   \mu_{st } ( d y ) \\
      \end{eqnarray*}
      On the other hand side,
      \begin{eqnarray*}
     \pi_c  \pi( t^{ l }   I_n ,  ( t^{\bf k }  , \mu_{st } ) ) f ( x ) &=& \lambda^l  \pi_c  \pi ( t^{\bf k }  , \mu_{st } )   f  \\
       &=& \lambda^l   \int_{ V_{ t^{\bf k } } }  f ( c^{-1} t^{ - \bf k } x +  y )   \mu_{st } ( d y ) .
      \end{eqnarray*}
      This proves the lemma.
       \hfill $\Box $

 The central extension splits in the subgroup $G$ given in (\ref{1.1}). It is easy to see that $G$ commutes with $ \widehat {SL_n} (F) $.
  For every character $\chi : G \to {\Bbb C}^*$, the space of $\chi$-coinvariants
  \begin{equation}
   \mathcal S(R^n_-)_{  \lambda , \chi } \deff \mathcal S(R^n_-)_{  \lambda } /  N (\chi )
   \end{equation}
   where $ N ( \chi ) $ is the linear span of elements $ \pi ( g ) f - \chi ( g ) f $,
  is a representation of $ \widehat {SL_n} (F) $.  Since
    $\pi_c$ commutes with $G$, it induces an isomorphism
    \[ \pi_c : \,  \mathcal S(R^n_-)_{  \lambda , \chi } \to \mathcal S(R^n_-)_{ |c|^n \lambda , \chi }. \]
 According to the analogy with the classical theta correspondence, $\mathcal S(R^n_-)_{  \lambda , \chi }$
  should be understood as a degenerate principal series representation. The principal series for  $\widehat {\frak sl}_2 $
   is considered in \cite{FZ} using a different method.

\

We will write an element  $x\in R_-^n = {\Bbb R}^n [ t^{-1} ] t^{-1} $ as
\[ x=  \sum_{ i=1}^{\infty }     x_{-i} t^{-i} ,\]
where $x_{-i}\in {\Bbb R}^n $. We write a function $f $ on $R_-^n $ as
\begin{equation}\label{2.14}   f (x_{-1} , x_{-2} , \dots  ) \end{equation}
accordingly. Sometimes it is more convenient to write each $x_{-i} \in {\Bbb R}^n $ as
$ x_{-i} = \sum_{ j =1}^n x_{-i}^j e_j $, where $e_j$ ($j=1, 2, \dots , n $) is the standard basis for ${\Bbb R}^n$, so
  \[  x=  \sum_{ i=1}^{\infty } \sum_{1 \leq  j \leq n }    x^{j}_{-i} e_j t^{-i} ,\]
$ f \in \mathcal S(R_-^n )$
is expressed as a function of variables  $x^{j}_{-i}$:
\begin{equation}\label{2.15}  f ( x^1_{-1} , x^{2}_{-1} , \cdots , x^n_{ -1}  ,   x^1_{-2} , x^{2}_{-2} , \cdots , x^n_{ -2}, \cdots ) . \end{equation}

\

\

\section{ Maximal Compact Subgroup of $\widehat {SL_n} ( F )$ and Its Fixed Gaussian Function.  }\label{section3}

In this section, we give a construction of a ``maximal compact subgroup'' $K$ of $\widehat {SL_n} (F )$
  and give a Gaussian function in  $ \mathcal S(R^n_-)_\lambda$ that is fixed by $ K $.

 We first define a ``maximal compact subgroup'' $O$ of $GL_n(F)$ by
 \begin{equation}
  O =\{ g(t) \in  {GL}_n (\mathbb R [ t , t^{-1} ]  )\; \vert \;  g(t)g(t^{-1})^T =I_n \},
 \end{equation}
 where $g(t^{-1})^T$ denotes the transpose of $g(t^{-1})$. Let $  d( t ) = {\rm det } \, g ( t) \in {\Bbb R} [ t , t^{-1} ]$. The condition
  $  g(t)g(t^{-1})^T = I_n$ implies that $ d ( t ) d ( t^{-1} ) = 1 $, so $ d ( t ) = \pm  t^k $ for some $ k \in {\Bbb Z}$. This proves
   $ O \subset GL_n (F)_0 $. We will prove the central extension $\widehat {GL_n} ( F)_0 $ splits over $ O$.

  Our definition of the ``orthogonal group" $O$ is justified by the following consideration.
We first define an inner product on $\mathbb R^n[t,t^{-1}]$. For any $u(t),v(t) \in \mathbb R^n[t,t^{-1}]$, define
    \begin{equation}\label{3.2}  (u(t),v(t)):=(u(t)^Tv(t^{-1}))_0     \end{equation}
   where $(w)_0$ denotes the constant term in $w$, for any $w\in \mathbb R [t,t^{-1}]$.
 It is clear that homogeneous pieces $ {\Bbb R}^n t^j$ are mutually orthogonal and
 the basis
 $\{e_i t^j\}_{1\leq i \leq n, j \in \mathbb Z}$ is orthonormal,   where $e_1,\cdots,e_n$ denotes the standard basis for $\mathbb R^n$.
 An element $ g ( t ) \in {GL}_n (\mathbb R [ t , t^{-1} ]  )$ preserves the inner product (\ref{3.2}) iff
 \[  (u(t)^Tv(t))_0 = (g(t)u(t),g(t)v(t))=(u(t)^Tg(t)^T g(t^{-1}) v(t^{-1}))_0\]
for all $ u ( t ) , v(t)   \in{\Bbb R}^n [ t, t^{-1} ] $
 iff $ g(t)^T g(t^{-1})= I_n $. Therefore $ O $ is the subgroup in $ GL_n  ( {\Bbb R} [ t , t^{-1}] ) $ that preserves the inner product
 (\ref{3.2}).

It is clear that $ g \in O $ implies that $ t^k g \in O$ for every $ k \in {\Bbb Z} $. So every $ g \in O$ can be
 written as $ g = t^k g' $ with $ g ' \in O \cap M_n ( {\Bbb R} [ t] )$.
For any $g \in O \cap M_n ( {\Bbb R} [ t] ) $,  we identify the quotient space ${V}_g$ with $ g^{-1}R^n \cap R^n_- \subset R^n_-$.
 The inner product on ${\Bbb R}^n [ t, t^{-1} ]$ restricts to an inner product on $V_g$. Let $\mu_g$ denote the Haar measure induced
 by this inner product in the sense that the volume of parallelotope spanned by an orthonormal basis in
 $V_g$ equals $1$.
From (\ref{2.4}), we can easily check
 \begin{equation}\label{3.3}
  \mu_{g_1}\ast \mu_{g_2}=\mu_{g_1g_2} \; \; \; {\rm for } \; g_1 ,g_2 \in O \cap M_n ( {\Bbb R} [ t] ). \end{equation}

For any $g\in O$, let $l$ be large enough so that $t^l g \in M_n( {\Bbb R} [ t])$, then define $s(g):=(t^{-l},(t^lg,\mu_{t^l g}))$.
 $s(g) $ is independent of the choice $l$ by (\ref{3.3}).
Hence $s(O)= \{ s ( g ) \, | \, g\in O \} $ is a subgroup of $\widehat{GL_n} (F )_0$ that is a lifting
 of $O \subset {GL}_n ( F )_0$.
 We define a
``maximal compact subgroup" $K$  for $\widehat {SL_n} ( F)$ as
\[   K =  \{ s ( g ) \, | \, g\in O , {\rm det} \, g =1 \}\]
 which is a lifting of $ O \cap SL_n (F )$.

\

Now we want to construct a nontrivial function in $ \mathcal S(R^n_-)_\lambda$ that is fixed by $K$.
 By simple observation we know for any function $ f \in\mathcal  S ( {\Bbb R}^n )$ with
\[  f ( 0 ) =1 , \; \; \;   \int_{{\Bbb R}^n} f ( x ) d x = \lambda \]
       then
  \[
  \phi_f ( x_{-1}^1 , \cdots , x_{-1}^n , \cdots ) :=
  \Pi_{m=1}^\infty  f ( x_{-m}^1 , \cdots , x_{-m}^n ) \in  \mathcal S(R^n_-)_\lambda
   \]
The infinite product makes sense since almost all of $x_{-m}^j$ are $0$ and $ f ( 0 ) =1 $.
We take
 \[  f_c ( x_1  , \cdots , x_n ) = e^{ - \pi  c  (x_1^2 + \cdots + x_n^2 ) }. \]
Then
\begin{equation}\label{3.4}
  \phi_{c}:=\phi_{f_c}= e^{ - \pi c \sum_{ j=1}^n \sum_{k=1}^\infty  {x_{  - k }^j }^2 } . \end{equation}
By calculation, we have
\[    \pi_t \phi_c ( x ) =    \lambda_c \phi_c ( x ) \]
where
\[    \lambda_c = \int  e^{ - \pi c (x_1^2 + \cdots + x_n^2 ) } dx_1 \cdots dx_n = c^{-\frac  n 2 } . \]
Fix $\lambda > 0 $, let $ c = \lambda^{ - \frac 2 n } $, then  $ \pi_t \phi_c ( x ) = \lambda \phi_c $.

\begin{prop} The function $\phi_c \in  \mathcal S(R^n_-)_\lambda $
   is fixed by $K$.
\end{prop}

\noindent {\it Proof.}
 Let $(t^{l},(g,\mu_{g}))\in K$, with $ g\in O \cap M_n(\mathbb R[t])$.
According to definition, we have
\[
\pi(t^{l},(g,\mu_g)) \phi_c(x)=\lambda^l \int _{V_g} e^{-\pi c((g^{-1} x)_-+y,  \, (g^{-1}x)_-+y )} \mu_g(dy),
\]
where we identify $V_g$ with $g^{-1}R^n\cap R^n_-$, and $(x)_- :=\pi_-(x)$.

 Since $g\in O \cap M_n(\mathbb R[t])$, obviously $g(t)^{-1}=g(t^{-1})^T \in M_n(\mathbb R[t^{-1}])$.
 This means $g(t)^{-1}x$ is indeed an element in $R^n_-$, and hence $(g^{-1}x)_-=g^{-1}x$.
 Then since $g$ preserves the $(\cdot,\cdot)$, we have
 \[
 ((g^{-1} x)_-+y,  \, (g^{-1}x)_-+y )=(x+gy,\, x+gy).
 \]
But $gy\in R^n$ as $y\in g^{-1}R^n \cap R^n_-$. The orthogonality of $x$ and $gy$ implies the right side of the above equation equals
\[
(x,\,x)+(y,\,y).
\]
Thus
\begin{align*}
\pi(t^{l},(g,\mu_g)) \phi_c(x) &=\lambda^l \int_{V_g} e^{- \pi c (y,y)}\mu_g(dy)   \, \phi_c(x)\\
&=\lambda ^l c^{-\frac{1}{2}\dim V_g} \, \phi_c(x).
\end{align*}
And as $\det t^{l}g=1$, the dimension of $V_g$ equals $-ln$. Combining with the fact that $\lambda=c^{-\frac{n}{2}}$, leads to the desired result.
 \hfill $\Box $

\

\

\section { Formula for action of  Affine Lie Algebra $\widehat{\frak {sl}_n} $.} \label{section4}
 In this section, we derive a formula for  the action of affine algebra of $\frak {sl}_n $ on $ {\cal S} ( R_{-}^n )_\lambda $ corresponding
 to the loop group action constructed in Section 2.  Let  $\mathfrak g =\frak {sl}_n  $  and $\hat {\mathfrak g} $ be the corresponding affine Kac-Moody algebra, so
 \[  \hat {\mathfrak g} = \frak {sl}_n \otimes {\Bbb C} [ t, t^{-1} ] + {\Bbb C} K \]
   with Lie bracket given by
   \[     [ a t^m , b t^n ] = [ a , b ] t^{ m+ n} + m \delta_{ m + n , 0 } ( a , b ) K \]
   where $ ( a , b ) ={\rm tr} \, a b $.
   By (\ref{2.5}),   $ GL_n ( R )$ acts on $  {\cal S} ( R_{-}^n )$,  its Lie algebra $ gl_n \otimes {\Bbb C} [ t ] $ also acts on the same space.
  We  write an function $f \in {\cal S} ( R_{-}^n )$ as (\ref{2.15}).  The action of $ gl_n \otimes {\Bbb C} [ t ] $ is given by first order differential operator
   as in the following lemma.

\begin{lemma}  Let $E_{uv} \in gl_n $ be the matrix with $(u,v)$-entry $1$ and other entries $0$.  For $j \geq 0 $,   $ E_{uv} t^j \in gl_n \otimes {\Bbb C} [ t ]$, let
\begin{equation}\label{4.1}
     \pi (  E_{uv} t^j )  =  - \sum_{ i = 1}^{\infty}   x^v_{ -i - j }   \partial_ {x^u_{ - i } }.
 \end{equation}
   This gives a representation of $  gl_n \otimes {\Bbb C} [ t ]$ on ${\cal S} ( R_-^n )$.
\end{lemma}

\noindent {\it Proof.}
 We first argue that $ \pi (  E_{uv} t^j ) f \in  {\cal S} ( R_{-}^n ) $ for $ f \in   {\cal S} ( R_{-}^n )$. It is
   enough to prove, for every $m \geq 1 $, the restriction of  $\pi (  E_{uv} t^j ) f $ to ${\Bbb R}^n t^{-1} + \dots +  {\Bbb R}^n t^{-m} $ is a Schwartz function.
    This restriction is
 \[      -\sum_{i=1}^{m-j} x^v_{-i-j} \partial_{x_{-i}^u} f(x_{-1} , x_{-2} , \dots , x_{-m} , 0 , \dots  ), \]
  which is in ${\cal S} [  {\Bbb R}^n t^{-1} + \dots +  {\Bbb R}^n t^{-m} ] $.
  It can be proved by a direct calculation that  $ [ \pi ( a ) , \pi ( b ) ]   = \pi ( [ a , b ] ) $ for $ a , b \in gl_n [ t ] $.
\hfill $\Box $

\

 The formula (\ref{4.1}) is derived from the formal calculation
 \[      \pi ( E_{ uv } t^j) f ( x ) =  \left( \frac d { d \epsilon}   \pi ( 1 + \epsilon  E_{ uv } t^j , \mu_{st} )  f ( x )  \right)|_{\epsilon = 0 }. \]

\[ (1+\epsilon E_{uv} t^j)^{-1}=I+(-\epsilon ) E_{uv }t^{-j}, \; \; \; {\rm mod} \; \epsilon^2  \]
\[  f (  ( 1 - \epsilon E_{uv} t^j ) x )
     = f ( x_{-1} -  \epsilon x_{ - 1 - j }^v e_u  ,  \cdots , x_{-m} - \epsilon  x_{ - m - j }^v   e_u,\cdots )\]
Taking derivative $\frac d { d \epsilon}|_{\epsilon = 0 } $, we derive (\ref{4.1}).

\

\begin{theorem} \label{theorem4.2}   Let $ \pi : \hat {\frak g} \to End ( {\cal S} ( R_-^n )_\lambda ) $ be the linear map given as follows:
$ \pi ( a t^j ) $ is {\rm (\ref{4.1})} for $ a\in \frak g $ and $ j \geq 0 $; for $ j > 0 $, $ u \ne v $,
 \begin{equation} \label{4.2}
    \pi(E_{uv}t^{-j})=-\lambda^{-j}\pi_t^j\sum_{i=1}^\infty x_{-i}^v \partial_{x_{-i-j}^u},
    \end{equation}
and
\begin{equation} \label{4.3}
     \pi ((E_{uu}-E_{vv})t^{-j})=-\lambda^{-j}\pi_t^j\sum_{i=1}^\infty (x_{-i}^u \partial_{x_{-i-j}^u} -x_{-i}^v \partial _{x_{-i-j}^v }) .
     \end{equation}
Then this gives  a representation of affine Kac-Moody algebra $\hat{\mathfrak g}$ on $ {\cal S} ( R_-^n )_\lambda  $ with central charge $1$.
\end{theorem}

\noindent {\it Proof.}  We have already proved in the proof of Lemma 4.1 that $ \pi ( a t^j ) $ for $j\geq 0 $ preserves the space of Schwartz functions.
The similar method applies to show that  $ \pi ( a t^{-j} ) $, $ j \geq 1 $
sends a Schwartz function on $ R_-^n$ to a Schwartz function on $ R_-^n$. Then
we need to prove that the operators $ \pi ( a t^j ) $ for $ a\in \frak g $ and $ j \geq 0 $ and operators in (\ref{4.2}), (\ref{4.3}) preserve the eigenspace
 $ {\cal S} ( R_-^n )_\lambda $. To this end,
 we need the relations:
\begin{equation}\label{4.4}
    x_{-i}^s  \pi_t = \pi_t  x_{-i-1}^s, \; \; \; \partial_{x_{-i}^s}  \pi_t = \pi_t \partial_{ x_{-i-1}^s} , \; \; \;  \pi_t^n  \partial_{x_{-i}^s}= 0  \; \; {\rm for}
   \; 1\leq i \leq n .
\end{equation}
Using
 \begin{equation}
   \pi_t f( x_{-1} , x_{-2} , \dots  ) =  \int_{ y \in {\Bbb R}^n } f ( y , x_{-1} , x_{-2} , \dots ) d y,
 \end{equation}
 one checks the first two relations directly. The third follows from the fact that
 \[    \int_{ \Bbb R} \partial_{x } f ( x ) d x = 0 \]
for a Schwartz function $ f \in {\cal S} ( {\Bbb R} )$.
  Using (\ref{4.4}), one can check the operators $\pi ( a t^i )$, $ j \in {\Bbb Z} $,  commute with $\pi_t$ therefore they preserves the eigenspace ${\cal S} ( R_-^n )_\lambda $.
  We verify
   $ \pi_t $ commutes with (\ref{4.2}).
  \begin{eqnarray*}
     \pi_t \pi(E_{uv}t^{-j}) &=& - \lambda^{-j}\pi_t^{j+1}\sum_{i=1}^\infty x_{-i}^v \partial_{x_{-i-j}^u}  \\
  &=& - \lambda^{-j} \pi_t^{j+1} x_{-1}^v   \partial_{x_{-1-j}^u}      -  \lambda^{-j}\pi_t^j\sum_{i=2}^\infty \pi_t    x_{-i}^v  \partial_{x_{-i-j}^u}\\
  &=& - \lambda^{-j} \pi_t^{j+1} \partial_{x_{-1-j}^u} x_{-1}^v    -  \lambda^{-j}\pi_t^j\sum_{i=2}^\infty    x_{-i+1}^v  \partial_{x_{-i+1-j}^u} \pi_t \\
  &=& -  \lambda^{-j}\pi_t^j\sum_{i=1}^\infty    x_{-i}^v  \partial_{x_{-i-j}^u} \pi_t \\
  &=&  \pi(E_{uv}t^{-j})    \pi_t
 \end{eqnarray*}
The other commutation relations can be verified similarly.
Next we prove $ [ \pi ( a ) , \pi ( b ) ] = \pi ( [ a , b ] ) $ for $ a, b \in \hat{\frak g } $ and the central charge $K$ acts as $1$.
 The case $ a , b \in \frak {sl}_n [ t ] $ has already appeared in Lemma 4.1.
 We compute the case $  a=   E_{vu} t $ and $ b = E_{uv} t^{-1} $. The other cases are similar. Notice that
 \begin{eqnarray*}
   \pi ( E_{u v } t^{-1} ) &=& -\lambda^{-1}\pi_t \sum_{i=1}^\infty x_{-i}^v \partial_{x_{-i-1}^u} \\
 &=& -\lambda^{-1}\pi_t^1 x_{-1}^v \partial_{x_{-2}^u} -  \lambda^{-1}\pi_t   \sum_{i=2}^\infty x_{-i}^v \partial_{x_{-i-1}^u} \\
 &=& -\lambda^{-1}\pi_t^1 x_{-1}^v \partial_{x_{-2}^u} -  \lambda^{-1}  \sum_{i=1}^\infty x_{-i}^v \partial_{x_{-i-1}^u}\pi_t  \\
  &=& -\lambda^{-1}\pi_t^1 x_{-1}^v \partial_{x_{-2}^u} -    \sum_{i=1}^\infty x_{-i}^v \partial_{x_{-i-1}^u}
  \end{eqnarray*}
  the last ``$=$" follows from the fact that $ \pi_t $ acts on $ {\cal S} ( R_-^n)_\lambda $ as $\lambda $.
  Now
   \begin{eqnarray}
 [ \pi ( E_{v u} t ) ,   \pi ( E_{u v } t^{-1} ) ] &=&  [  - \sum_{ i = 1}^{\infty}   x^u_{ -i - 1 }   \partial_ {x^v_{ - i } }  ,
- \lambda^{-1} \pi_t  x^v_{-1}  \partial_ {x^u_{ - 2 } } - \sum_{ i = 1 }^{\infty}   x^v_{ -i  }   \partial_ {x^u_{ - i-1 } }]\\
   &=& \sum_{ i = 1}^{\infty} [   x^u_{ -i - 1 }   \partial_ {x^v_{ - i } }  ,   \lambda^{-1} \pi_t    x^v_{ -1   }   \partial_ {x^u_{ -2 } } ]
     + \sum_{ i =1}^\infty [ x^u_{ -i - 1 }   \partial_ {x^v_{ - i } } ,  x^v_{ -i  }   \partial_ {x^u_{ - i-1 } }]   \nonumber
    \end{eqnarray}

 \begin{equation}
 \sum_{ i =1}^\infty [ x^u_{ -i - 1 }   \partial_ {x^v_{ - i } } ,  x^v_{ -i  }   \partial_ {x^u_{ - i-1 } }]
  = \sum_{ i = 1}^{\infty}   x^u_{ -i-1  }   \partial_ {x^u_{ - i-1  } } -  \sum_{ i = 1 }^{\infty}   x^v_{ -i  }   \partial_ {x^v_{ - i } }.
  \end{equation}

    \begin{eqnarray}
 &&  \sum_{ i = 1}^{\infty} [   x^u_{ -i - 1 }   \partial_ {x^v_{ - i } }  ,   \lambda^{-1} \pi_t    x^v_{ -1   }   \partial_ {x^u_{ -2 } } ] \nonumber \\
 &&=   \sum_{ i = 1}^{\infty}x^u_{ -i - 1 }   \partial_ {x^v_{ - i } }   \lambda^{-1} \pi_t    x^v_{ -1   }   \partial_ {x^u_{ -2 } }
      - \sum_{ i = 1}^{\infty} \lambda^{-1} \pi_t    x^v_{ -1   }   \partial_ {x^u_{ -2 } }   x^u_{ -i - 1 }   \partial_ {x^v_{ - i } } \nonumber \\
&&=   \sum_{ i = 1}^{\infty} \lambda^{-1} \pi_t   x^u_{ -i - 2 }   \partial_ {x^v_{ - i -1  } }     x^v_{ -1   }   \partial_ {x^u_{ -2 } }
      - \sum_{ i = 1}^{\infty} \lambda^{-1} \pi_t    x^v_{ -1   }   \partial_ {x^u_{ -2 } }   x^u_{ -i - 1 }   \partial_ {x^v_{ - i } } \nonumber \\
 && =  - \lambda^{-1} \pi_t    x^v_{ -1   }   \partial_ {x^u_{ -2 } }   x^u_{ - 2 }   \partial_ {x^v_{ - 1 } }   \nonumber \\
 && =  - \lambda^{-1} \pi_t     (  \partial_ {x^v_{ - 1 } } x^v_{ - 1 } -1 )  (  x^u_{ - 2 }    \partial_ {x^u_{ -2 } } + 1 )     \nonumber \\
  && =\lambda^{-1} \pi_t  (  x^u_{ - 2 }    \partial_ {x^u_{ -2 } } + 1 )    \nonumber \\
  &&=  x^u_{ - 1 }    \partial_ {x^u_{ -1 } } + 1
 \end{eqnarray}
   Combine (4.6), (4.7) and (4.8), we get
  \[  [ \pi ( E_{v u} t ) ,   \pi ( E_{u v } t^{-1} ) ] =  \pi ( E_{vv} - E_{uu} ) + 1 \]
   This proves $K=1$.
   \hfill $\Box $

\

Next we explain how (\ref{4.2}) is derived from the group action in Section 2. The case (\ref{4.3}) is similar.
 To ease the notations, we assume $ n=2$ and $ E_{uv}=E_{12}$. Notice that
 \[         E_{12} t^{-j} = \left( \begin{matrix} t^{-j} & 0 \\ 0 & t^j \end{matrix} \right)
                      E_{12} t^j  \left( \begin{matrix} t^{j} & 0 \\ 0 & t^{-j} \end{matrix} \right)
                  =  g \, E_{12} t^j \, g^{-1} \]
 and $ g \in O $.  We show $  \pi ( E_{12} t^{-j} )$ given in (\ref{4.2}) is equal to
  $ \pi ( g )  \pi (  E_{12} t^j ) \pi ( g^{-1} ) $,
   where $ \pi (  E_{12} t^j ) $ is given by (\ref{4.1}).
   \[ \pi ( g ) = \pi (   t^{-j} , ( t^{j} g , \mu_{ t^{j} g  } ) ),
    \, \, \;   \pi ( g^{-1}  )= \pi (   t^{-j} , ( t^{j} g^{-1} , \mu_{ t^{j} g^{-1}  } ) )\]
 To write  $\pi ( g )$ and $ \pi ( g^{-1}  )$ more explicitly, we introduce operators
 $ \pi_{ t, 1} $ and $ \pi_{ t, 2 }$ on ${\cal S} ( R_-^2 ) $.
 \[ \pi_{ t, 1} f( x^1_{-1} , x^2_{-1} , x^1_{-2} , x^2_{-2} ,\dots )
   = \int_{\Bbb R} f( y , x^2_{-1} , x^1_{-1} , x^2_{-2} , x^1_{-2} , x^2_{-3} , \dots ) dy \]
   \[ \pi_{ t, 2} f( x^1_{-1} , x^2_{-1} , x^1_{-2} , x^2_{-2} ,\dots )
   = \int_{\Bbb R} f( x^1_{-1} , y ,  x^1_{-2} ,    x^2_{-1} , x^1_{-3} , x^2_{-2}  \dots ) dy \]
 Then
\[  \pi ( g ) =\lambda^{-j}  \pi_{ t, 2}^{2j} , \; \;  , \pi ( g^{-1} ) =\lambda^{-j}  \pi_{ t, 1}^{2j}.\]
 Notice that the relations  $\pi_t = \pi_{ t , 1 } \pi_{ t , 2 } = \pi_{ t , 2 }\pi_{ t , 1 } $.
 \begin{equation}
    x_{-i}^u  \pi_{t , v} = \pi_{t , v} x_{-i}^u ,  \; \; \partial_{x_{-i}^u} \pi_{t , v} = \pi_{t , v} \partial_{x_{-i}^u} \; \; {\rm when}
      \; u \ne v .
    \end{equation}
   and
  \begin{equation}
  x_{-i}^s \pi_{t , s}  = \pi_{t , s}  x_{-i-1}^s, \; \; \; \partial_{x_{-i}^s}  \pi_{t , s}
  = \pi_{t , s} \partial_{ x_{-i-1}^s} , \; \; \;  \pi_{t, s}^n  \partial_{x_{-i}^s}= 0  \; \; {\rm for}
   \; 1\leq i \leq n .
\end{equation}
\[ \pi ( g )  \pi (  E_{12} t^j ) \pi ( g^{-1} )
=  \lambda^{-2j}  \pi_{ t, 2}^{2j}  \left( -\sum_{i=1}^{\infty} x^2_{-i-j} \partial_{x_{-i}^1} \right)  \pi_{ t, 1}^{2j}\]
We apply relations to move $ \pi_{ t, 2}^j $ to the right and $ \pi_{ t, 1}^{j}$ to the left,
we find the result is (\ref{4.2}).

\

Parallel to the result in the section 2, we will show the function $\phi_c$ is killed by a subalgebra of $\hat{\mathfrak g}$.

\begin{prop} \label{proposition 4.3}
The function   $\phi_c \in \mathcal S ( R^n_- )_\lambda $ as in (\ref{3.4})
is killed by the fixed points of Chevalley involution of $\hat{\mathfrak g}$.
\end{prop}

\noindent{\it Proof.} The fixed points of Chevalley involution $w$ is
${\rm Span}\{ at^m -a^T t^{-m}\; \vert\; a\in \frak {sl}_n,  m\in \mathbb Z\}$,
where $w(a t^m )=-a ^T t^{-m},\; w(K)=-K$. So suffices to prove $h_c$
is killed by ${\rm Span}\{ at^m -a^T t^{-m}\; \vert\; a\in \frak {sl}_n,  m\in \mathbb Z\}$.

We first verify $\pi(E_{uv} t^m -E_{vu} t^{-m}) \phi_c=0$ for $u\neq v$.
\begin{align*}
\begin{split}
 \pi(E_{uv}t^{m} ) \phi_c  =& -\sum_{i=1}^\infty x^{v}_{-i-m} \partial _{x^{u}_{-i}} \phi_c\\x
 =& 2\pi \lambda^{-\frac{2}{n}}\sum_{i=1}^\infty  x^{v}_{-i-m} x^{u}_{-i} \phi_c\\
 \pi(E_{vu}t^{-m}) \phi_c  =& -\lambda^{-m} \pi_t^m \sum_{i=1}^\infty x^{u}_{-i} \partial_{x^{v}_{-i-m}} \phi_c\\
=&-\lambda^{-m} \int_{V_{t^m}}\Big(\sum_{i=1}^m y^{u}_{m-i} (-2\pi \lambda^{-\frac{2}{n}}) x^{v}_{-i}+
 \sum_{i=m+1}^\infty x^{u}_{m-i} (-2\pi \lambda^{-\frac{2}{n}} )x^{v}_{-i} \Big) \\
 &\cdot e^{-\pi \lambda^{-\frac{2}{n}}  ((y^{1}_{m-1})^2 +\cdots +(y^{n}_{m-1})^2 +\cdots (y^{1}_{0})^2+\cdots +(y^{n}_{0})^2+(x^{1}_{-1})^2+\cdots) }
 d y^{1}_{m-1} \cdots  dy^{n}_{0}   \\
 =&2\pi \lambda^{-\frac{2}{n} -m }  \sum_{i=1}^m x^v_{-i}\phi_c \int   y^{u}_{m-i} e^{-\pi \lambda^{-\frac{2}{n}}((y^{1}_{m-1})^2 +\cdots +(y^n_{0})^2)}
 d y^1_{m-1} \cdots  dy^n_0 \\
 & +2\pi \lambda^{-\frac{2}{n}-m} \int e^{-\pi \lambda^{-\frac{2}{n}} ((y^1_{m-1})^2 +\cdots +(y^n_0)^2) } d y^1_{m-1}\cdots d y^n_0\cdot
 \sum_{i=m+1}^\infty x^u_{m-i} x^v_{-i} \phi_c\\
 =& 2\pi \lambda^{-\frac{2}{n}-m}   (\int e^{-\pi \lambda^{-\frac{2}{n}}  y^2} dy )^{mn} \sum_{i=1}^\infty x^u_{-i}x^v_{-i-m} \phi_c\\
=&2\pi \lambda ^{-\frac{2}{n}} \sum_{i=1}^\infty x^u_{-i}x^v_{-i-m}\phi_c
\end{split}
\end{align*}

Hence
\[
\pi(E_{uv}t^m-E_{vu}t^{-m} ) \phi_c=0.
\]

For the case $\pi((E_{uu}-E_{vv}) (t^m-t^{-m})) \phi_c=0$, similar calculation shows
\begin{align*}
\begin{split}
\pi((E_{uu}-E_{vv})t^m ) \phi_c
=& 2\pi \lambda^{-\frac{2}{n}} \sum_{i=1}^\infty (x^u_{-i-m} x^u_{-i} -x^v_{-i-m} x^v_{-i} ) \phi_c\\
=&\pi((E_{uu}-E_{vv})t^{-m}) \phi_c
\end{split}
\end{align*}

Hence
\[
\pi((E_{uu}-E_{vv}) (t^m-t^{-m})) \phi_c=0.
\]
\qed

\

\

\

\section{Relations With Highest Weight Modules}\label{section5}
In this section, we show that certain highest weight modules of $\hat {\frak g}$ appear in the dual representation of
the representation  ${\cal S} ( R^n_-)_\lambda $.

Let  $\mathcal S(R^n_-)_\lambda ^\ast$ be the full dual space of $\mathcal S(R^n_-)_\lambda$.
It is representation of affine Lie algebra $\hat {\frak g}$ under the action
\begin{equation} \label{5.1}
\la   \pi^\ast (a) \varphi ,f \ra := -  \la \varphi , \pi(a) f \ra.
\end{equation}
It is clear that the level of $\mathcal S(R^n_-)_\lambda ^\ast$ is $-1$ as $K$ acts as $-1$ by (\ref{5.1}).
This representation is too big to be interesting. However we show that some natural linear functionals are highest weight vectors.

Let $\mathfrak h$ be its Cartan subalgebra of $\frak g = \frak {sl}_n$
consisting of diagonal matrices, $n^+$ ($n^-$) be the strictly upper (lower) triangular matrices, so we have triangular decomposition
$ \frak g = n^- \oplus \mathfrak h \oplus n^+$. The corresponding triangular decomposition for $\hat {\frak g} $ (see \cite{Ka}) is
\[    \hat {\frak g} =  \hat n^- \oplus \hat {\frak h}  \oplus \hat n^{ + } \]
with
 \[    \hat {\mathfrak h}=\mathfrak h+ {\Bbb C} K ,   \; \; \;  \hat n^- = n_- \oplus \frak g [ t^{-1} ] t^{-1} ,  \; \; \;  \hat n^+ = n_+ \oplus \frak g [ t ] t . \]

We define linear functional $I_k$ for $1\leq k \leq n$ as follows, for  $f\in \mathcal S(R^n_-)_\lambda$,
\begin{equation} \label{5.2}
\la I_k , f \ra =  I_k(f) := \int_{\mathbb R^k } f(x_{-1}^1,\cdots,x_{-1}^k ,0,0\cdots ) dx_{-1}^1\cdots dx_{-1}^k
 \end{equation}
where $d$ is the Lebesgue measure on $\mathbb R$.

We consider $I_n$ first. For any $E_{uv} t^j, \; j\geq 1$,
\begin{equation}\label{5.3}
\la   \pi^\ast (E_{uv}t^j) I_n,f \ra = -  \la I_n , \pi(E_{uv}t^j) f \ra
                                = \la I_n ,\sum_{i=1}^\infty x_{-i-j}^v \partial_{x_{-i}^u} f \ra
  \end{equation}
 by (\ref{4.1}), which is $0$ by (\ref{5.2}). Hence $\pi^\ast (E_{uv} t^j) I_n=0,$ for $j\geq 1$.
Similarly by (\ref{4.1}),  we have
\begin{equation}\label{5.4}
\la \pi^\ast (E_{uv}) I_n,f \ra = \int _{\mathbb R^n} x_{-1}^v \partial _{x_{-1}^u }f(x_{-1},0,0\cdots) d x_{-1}
\end{equation}
The above formula implies immediately $\pi^\ast (E_{uv}) I_n=0$ for $u\neq v$ since $f$ is a Schwartz function.
Also it implies $\pi^\ast (E_{uu})I_n= -I_n$ due to $x_{-1}^u \partial_{x_{-1}^u}= -[\partial_{x_{-1}^u},x_{-1}^u]+
\partial_{x_{-1}^u} x_{-1}^u=-1+\partial_{x_{-1}^u} x_{-1}^u$.
Hence  $\pi^\ast (E_{uu}-E_{vv})$ also kills $I_n$.
Thus we have shown the dual action of any element in $\frak {sl}_n(R)=\frak {sl}_n(\Bbb R[[t]]) $ kills $I_n$.
Moreover, the central charge $K$ acts as $-1$ as mentioned earlier.

Similarly one can show $\pi^\ast (E_{uv}t^j)$ for $j\geq 1$ kills $I_k$.
Now we consider $\pi^* (  E_{uv} )$ for $ u \leq v $. If $v >k$, we see that $\pi^\ast(E_{uv}) I_k=0$ immediately from the definition of $I_k$.
For $u,v\leq k$, similar to (\ref{5.4}), we derive that $\pi^\ast(E_{uv})$ kills $I_k$ when $u\neq v$ and acts as $-1$ when $u=v$.
To conclude, we obtain the following theorem

\begin{theorem} \label{theorem 5.1}
For $k=1,\cdots,n$, $I_k \in   \mathcal S(R^n_-)_\lambda ^\ast  $
is a highest weight vector with highest weight $\Lambda_k $  given as follows:
\begin{align*}
\la \Lambda_k,E_{ii}-E_{i+1,i+1}\ra  &=0, \;\;  \;\; \;\text{ for } i\neq k, \\
 \la \Lambda_k, E_{kk}-E_{k+1,k+1} \ra  &=-1 ,\;\; \text{ for } 1\leq k\leq n-1 \\
 \la \Lambda_k,  K \ra &=-1.
 \end{align*}
\qed
\end{theorem}

\

\section{Whittaker Functionals  }

Let's denote the Lie algebra of linear map on $n\times n$ matrix space  $ M_n  $ by $ gl_{ n^2}$
and the group of invertible linear transformations on $M_n $ by $ GL_{n^2}$.
We have the Lie algebra embedding  $ \frak {sl}_n \oplus \frak {sl}_n \to   gl_{ n^2} $ by the action
\[    ( a , b ) x = a x + x b^t \]
and the corresponding  group embedding
 $ SL_n \oplus SL_n \to   GL_{ n^2} $ is given  by the action
  \[  ( g_1 , g_2 ) x = g_1 x g_2^t . \]
 As in Section 2 and 4,  we know
 $ {\cal S} ( M_n ( {\Bbb R} ) [ t ^{-1} ] t^{-1} )_\lambda  $ is a representation of $ \widehat{ SL_{n^2}} ( {\Bbb R}) $ and
  affine Lie algebra $ \widehat {\frak {sl}_n } $, therefore is a representation of
   $  \widehat{ SL_{n}} ({\Bbb R}) \times   \widehat{ SL_{n}} ( {\Bbb R}) $ and Lie algebra
       $  \widehat{ \frak {sl}_{n}} \oplus  \widehat{ \frak {sl}_{n}} $ by restriction, its level is $ n $.
    The dual space ${\cal S} ( M_n ( {\Bbb R} ) [ t ^{-1} ] t^{-1} )_\lambda^*  $,
    as a representation of    $  \widehat{ \frak {sl}_{n}} \oplus  \widehat{ \frak {sl}_{n}} $, has level $-n$, which is the critical level as the dual Coxeter
   number of $ \frak {sl}_n$ is $n $.  We will construct Whittaker functionals  ${\cal S} ( M_n ( {\Bbb R} ) [ t ^{-1} ] t^{-1} )_\lambda^*  $.
    The Whittaker functionals for $\widehat{sl}_2$ are studied abstractly in \cite{ALZ}, our construction provides
     a concrete realization.

  We denote the action of the first (second) copy of $\widehat{\frak {sl}_n} $ in $ \widehat{\frak {sl}_n} \oplus \widehat{\frak {sl}_n}$ by $ \pi_1 (a ) $ ($\pi_2 ( a ) $)
  and the similar notations $\pi_1 ( g) $ and $ \pi_2 ( g)$ will be used for loop groups.
  Same convention will be used for the corresponding dual action $\pi_1^\ast$ and $\pi_2^\ast$ on the dual space
  ${\cal S} ( M_n ( {\Bbb R} ) [ t ^{-1} ] t^{-1} )_\lambda^*$.

 As in Section 2, we will write an element  $x\in M_n(\Bbb R) [t^{-1}]t^{-1} $ as
\[ x=  \sum_{ i=1}^{\infty }     x_{-i} t^{-i} ,\]
where $x_{-i}\in M_n(\Bbb R) $.
We write a function $f $ on $ M_n(\Bbb R) [t^{-1}]t^{-1} $ as
\begin{equation}\label{6.1}   f (x_{-1} , x_{-2} , \dots  ) \end{equation}
or as
 a function of variables  $x_{-i}^{u,v}$:
\begin{equation}\label{6.2}
f ( x^{1,1}_{-1} , x_{-1}^{1,2} , \cdots , x_{-1}^{ n,n}  ,   x_{-2}^{1,1} , x_{-2}^{1,2} , \cdots , x_{-2}^{ n,n}, \cdots ) .
\end{equation}
 where $x_{-i}^{u,v}$ denotes $(u,v)$-entry of $x_{-i} $.

  By (4.1) (4.2) (4.3), we have the following formulas
\begin{equation}
\pi_1 (E_{uv}t^j)=-\sum_{i=1}^\infty \sum_{s=1}^n x^{v,s}_{-i-j}\partial _{x^{u,s}_{-i}},
\end{equation}
and for $j>0,u\neq v$,
\begin{equation}
\pi_1 (E_{uv}t^{-j})=-\lambda^{-j} \pi_t^j \sum_{i=1}^\infty \sum_{s=1}^n x^{v,s}_{-i} \partial_{x^{u,s}_{-i-j}},
\end{equation}
and
\begin{equation}
\pi_1 ((E_{uu}-E_{vv})t^{-j})=-\lambda^{-j}\pi_t^j \sum_{i=1}^\infty \sum_{s=1}^n(x^{u,s}_{-i}\partial_{x^{u,s}_{-i-j}}-x^{v,s}_{-i} \partial_{x^{v,s}_{-i-j}} )
\end{equation}
\begin{equation}
\pi_2 (E_{uv}t^j)=-\sum_{i=1}^\infty \sum_{s=1}^n x^{s,v}_{-i-j}\partial _{x^{s,u}_{-i}},
\end{equation}
and for $j>0,u\neq v$,
\begin{equation}
\pi_2 (E_{uv}t^{-j})=-\lambda^{-j} \pi_t^j \sum_{i=1}^\infty \sum_{s=1}^n x^{s,v}_{-i} \partial_{x^{s,u}_{-i-j}},
\end{equation}
and
\begin{equation}
\pi_2 ((E_{uu}-E_{vv})t^{-j})=-\lambda^{-j}\pi_t^j \sum_{i=1}^\infty \sum_{s=1}^n(x^{s,u}_{-i}\partial_{x^{s,u}_{-i-j}}-x^{s,v}_{-i} \partial_{x^{s,v}_{-i-j}} )
\end{equation}

In this setting, take $\phi_c=e^{-\pi c\sum_{1\leq u,v\leq n} \sum_{k=1}^\infty (x^{u,v}_{-k})^2 }$ and $c=\lambda ^{-\frac{2}{n^2}}$.
We can easily show that $\phi_c$ is fixed by the maximal compact subgroup
 $ K \times K $ and killed by its Lie algebra as in Section 3.

Next we introduce Whittaker functionals in  ${\cal S} ( M_n ( {\Bbb R} ) [ t ^{-1} ] t^{-1} )_\lambda^*  $.
Our formula is motivated by the following finite dimensional case.  For simplicity we consider the case $ GL_3 ( {\Bbb R} ) \times GL_3 ( {\Bbb R} ) $.
It acts on the Schwartz space  $ {\cal S} ( M_3 ( {\Bbb R} ))$ as follows,
\begin{equation} \label{6.7}
  \pi ( g_1 , g_2 ) f ( x ) =   f ( g_1^{-1} x g_2^{-T} ) . \end{equation}

We still denote the action of the first (second) copy of $ GL_3 ( {\Bbb R} ) $ in
 $GL_3 ( {\Bbb R} ) \times GL_3 ( {\Bbb R} )$ by $ \pi_1 ( g) $ ($\pi_2 ( g ) $).

Let $ B$ be the standard Borel subgroup of $GL_3 ( {\Bbb R})$, i.e., $B$ consists of upper triangular matrices.
 Then $ B^T$, the lower triangular matrices, is an opposite Borel subgroup. $ B_2 = B \times B^T $ is a Borel subgroup of
 $GL_3 ( {\Bbb R} ) \times GL_3 ( {\Bbb R} ) $.
For $c = (c_1 , c_2) \in {\Bbb R}^2  $, we define a linear functional
 $\Phi _c: \cal S(M_3(\mathbb R)) \to   {\Bbb C}$ by
\begin{equation} \label{6.8}
\Phi_c(f):= \int_{\mathbb R^8}
                   f   \left(   \begin{pmatrix} x_{11} & x_{12} & x_{13}\\ x_{21} & x_{22} & x_{23} \\0 & x_{32} &x_{33}\end{pmatrix} \right)
                   e^{2\pi i (c_1 \frac{x_{11}}{x_{21}}+c_2 \frac{x_{22}} {x_{32}})} dx
\end{equation}
Let $\chi_c$ be the character of the  unipotent radical $U$ of $B$ defined by
 \[ \chi _c (u):= e^{-2\pi i(c_1u_1+c_2u_2)}, \;\;\;\; \text{ for }
      u= \begin{pmatrix}
              1 & u_1 &w \\
              0 & 1 & u_2\\
              0 & 0& 1
           \end{pmatrix}.
 \]
For such $u\in U$, we have
\begin{eqnarray}\label{6.9}
uX & =& \begin{pmatrix}
                    1 & u_1 &w \\
                    0 & 1 & u_2\\
                    0 & 0& 1
               \end{pmatrix}
               \begin{pmatrix}
                    x_{11} & x_{12} & x_{13}\\
                    x_{21} & x_{22} & x_{23} \\
                    x_{31} & x_{32} &x_{33}
               \end{pmatrix} \nonumber  \\
     & = &   \begin{pmatrix}
                    x_{11}+u_1x_{21}+wx_{31} & x_{12} +u_1x_{22}+w x_{32}& x_{13}+u_1x_{23}+wx_{33}\\
                    x_{21}+u_2x_{31} & x_{22} +u_2 x_{32}& x_{23}+u_2x_{33} \\x_{31} & x_{32} &x_{33}
                \end{pmatrix}
\end{eqnarray}
We check that
 \begin{equation} \label{6.10}
(\pi_1^\ast(u)\Phi_c,  f)=\chi_c(u)(\Phi_c,f),
\end{equation}
where $\pi_1^\ast$ is the dual action of $\pi_1$ on the full dual space of $\cal S(M_3(\Bbb R))$.
Indeed, by applying the functional $\Phi_c$ to $\pi_1(u)^{-1} f$ and using (\ref{6.9}), the left side equals
\[
\int_{\mathbb R^8} f
     \left(   \begin{pmatrix} x_{11}+u_1x_{21} & x_{12} +u_1x_{22}+w x_{32}& x_{13}+u_1x_{23}+wx_{33}\\
                                        x_{21} & x_{22} +u_2 x_{32}& x_{23}+u_2x_{33} \\
                                        0 & x_{32} &x_{33}
               \end{pmatrix}
     \right)   e^{2\pi i (c_1 \frac{x{11}}{x_{21}}+c_2 \frac{x_{22}} {x_{32}})} dx
  \]
 After changing variables, the above formula will be turned into
\[
\int_{\mathbb R^8}
     f  \left(
               \begin{pmatrix}
                    x_{11} & x_{12} & x_{13}\\ x_{21} & x_{22} & x_{23} \\0 & x_{32} &x_{33}
               \end{pmatrix}
        \right)
     e^{2\pi i (c_1 \frac{x_{11}-u_{1}x_{21}}{x_{21}}+c_2 \frac{x_{22}-u_2x_{32}} {x_{32}})} dx,
\]
which is indeed the right side of (\ref{6.10}).

Thus the functional $\Phi_c $ in (\ref{6.8}) is a Whittaker functional for the first copy of $GL_3(\Bbb R)$.
Similarly we will also give the Whittaker functional for the second copy of $GL_3(\Bbb R)$ as follows:
For $d = ( d_1,d_2) = \Bbb R^\ast$, define a linear functional $\Phi_d': \cal S(M_3(\Bbb R)) \to \Bbb C$ by
\begin{equation} \label{6.11}
\Phi_d' (f):= \int_{\mathbb R^8}
                   f   \left(   \begin{pmatrix} x_{11} & x_{12} & x_{13}\\ x_{21} & x_{22} & x_{23} \\0 & x_{32} &x_{33}\end{pmatrix} \right)
                   e^{2\pi i (d_1 \frac{x{22}}{x_{21}}+d_2 \frac{x_{33}} {x_{32}})} dx
\end{equation}
Let $\chi_d$ be a character of the  unipotent radical $U^T$ in $B^T$ defined by
 \[ \chi _d(v):= e^{-2\pi i(d_1v_1+d_2v_2)}, \;\;\;\; \text{ for }
      v= \begin{pmatrix}
              1 & 0 &0 \\
               v_1& 1 &0\\
               z & v_2& 1
           \end{pmatrix}
 .\]
In this case we have
 \begin{equation} \label{6.12}
(\pi_2^\ast(v)\Phi_d',  f)=\chi_d (v)(\Phi_d',f),
\end{equation}
where $\pi_2^\ast$ is the dual action of $\pi_2$ on the full dual space of $\cal S(M_3(\Bbb R))$.

We consider the loop group case $\widehat{SL_n}(F) \times \widehat{SL_n}(F)$.
It acts on the space $\cal S(M_n(\Bbb R)[t^{-1}]t^{-1})_\lambda$ with the action explained in the beginning of section 6.
Let $\hat B$ be the Borel subgroup consisting of elements $b\in \widehat{SL_n}(R)$ such that $b(0)$ is upper triangular.
Then $\hat B^T$, is an opposite Borel subgroup. And $\hat B\times \hat B^T$ is a Borel subgroup of
$\widehat{SL_n}(F) \times \widehat{SL_n}(F)$.

For $c = ( c_1,c_2,\cdots,c_n ) \in \Bbb R^n $, we define a linear functional $\Psi_c :\cal S(M_n(\mathbb R[t^{-1}]t^{-1}))\to \Bbb C$ by
\begin{align}\label{6.13}
\Psi_c(f)=& \int_{ {\Bbb R}^{n^2}} f
\left(
  \begin{pmatrix}
      x_{-1}^{1,1} & \cdots &x_{-1}^{1,n-1}&x_{-1}^{1,n}\\
      x_{-1}^{2,1} & \cdots & x_{-1}^{2,n-1} & x_{-1}^{2,n} \\
      & \ddots   & \vdots & \vdots\\ 0 && x_{-1}^{n,n-1} & x_{-1}^{n,n}
  \end{pmatrix} ,
  \begin{pmatrix}
      0&\cdots & 0 & x_{-2}^{1,n}\\
      0 & \cdots & 0 & 0\\
      \vdots & \ddots & \vdots & \vdots\\
      0 & \cdots & 0 & 0
  \end{pmatrix},
  0,\cdots
\right)   \nonumber \\
& \cdot e^{2\pi i \left(c_1 \frac{x_{-1}^{1,1}}{x_{-1}^{2,1}} +\cdots +c_{n-1}\frac{x_{-1}^{n-1,n-1}}{x_{-1}^{n,n-1}} +
                   c_n \frac{x_{-1}^{n,n}}{x_{-2}^{1,n}}
\right)} dx
\end{align}
Let $\hat \chi_c $ be a character on the unipotent subgroup of $B$ defined as follows:
For any $u\in B$ such that $u(0)$ is upper triangular unipotent, that is
\[
     u=   \begin{pmatrix} 1 & u_1 & \cdots&\ast \\  & \ddots & \ddots &\vdots \\  &  & \ddots & u_{n-1}\\ && &1\end{pmatrix} +
           \begin{pmatrix}
              v_1&\cdots &\cdots\\
              v_2 &\cdots&\cdots\\
              \vdots&\ddots&\ddots\\
              v_n &\cdots &\cdots
            \end{pmatrix}t+\cdots,
  \]
 then define
\[
  \hat \chi_c (u):= e^{-2\pi i (c_1u_1+\cdots + c_{n-1}u_{n-1} +c_nv_n)},
\]
where $u_1,\cdots,u_{n-1}$ are superdiagonal elements of $u\in B$, and $v_n$ is the $(n,1)$-entry of $u'(0)$.

Now we show that $\Psi_c$ is indeed a Whittaker functional of $\widehat{SL_n}(F)$
 under the group action $\pi_1$ on $\cal S(M_n(\Bbb R)[t^{-1}]t^{-1})_\lambda$ in the sense that
\begin{equation}\label{6.14}
(\pi_1^\ast(u)\Psi_c, f)=\hat\chi_c(u)(\Psi_c,f).
\end{equation}
Indeed, by $(\pi_1^\ast(u)\Psi_c, f)=(\Psi_c, \pi_1(u)^{-1}f)$,  the left side equals
 \[
 \begin{split}
    \int f
           \left(
                     \begin{pmatrix}
                            x_{-1}^{1,1}+u_1 x_{-1}^{2,1} & \cdots & x_{-1}^{1,n-1}+u_1x_{-1}^{2,n-1}+\cdots & x_{-1}^{1,n}
                            +u_1x_{-1}^{2,n}+\cdots+v_1x_{-2}^{1,n} \\
                            x_{-1}^{2,1} & \cdots &x_{-1}^{2,n-1}+u_2x_{-1}^{3,n-1}+\cdots & x_{-1}^{2,n}
                            +u_2x_{-1}^{3,n}+\cdots +v_2 x_{-2}^{1,n}\\
                           & \ddots &\vdots &\vdots\\ 0 &  & x_{-1}^{n,n-1} & x_{-1}^{n,n}+v_nx_{-2}^{1,n}
                      \end{pmatrix} , \right.\\ \left.
                     \begin{pmatrix}
                               0 &\cdots &0 & x_{-2}^{1,n} \\ 0 &\cdots & 0 & 0\\ \vdots & \ddots & \vdots &\vdots \\ 0 &\cdots & 0 & 0
                       \end{pmatrix} ,   0,\cdots
             \right)
        \cdot e^{2\pi i \left( c_1 \frac{x_{-1}^{1,1}}{x_{-1}^{2,1}}+
        \cdots c_{n-1} \frac{x_{-1}^{n-1,n-1}}{x_{-1}^{n,n-1}}+c_n\frac{x_{-1}^{n,n}}{x_{-2}^{1,n}} \right) }dx
 \end{split}
\]
 Since the Jacobi matrix is upper triangular with all of diagonal entries equals one, hence the Jacobi determinant equals $1$.
 And by changing the variables, the above formula will be turned into
\begin{multline*}
  \int f
             \left(
                 \begin{pmatrix}
                          x_{-1}^{1,1} & \cdots &x_{-1}^{1,n-1}  &x_{-1}^{1,n}\\
                          x_{-1}^{2,1} & \cdots & x_{-1}^{2,n-1} & x_{-1}^{2,n} \\
                                              & \ddots & \vdots             & \vdots\\
                          0                  &            & x_{-1}^{n,n-1} & x_{-1}^{n,n}
                 \end{pmatrix} ,
                \begin{pmatrix}
                         0          &\cdots &          0 & x_{-2}^{1,n}\\
                         0          & \cdots &         0 &   0\\
                          \vdots & \ddots & \vdots & \vdots\\
                         0          & \cdots &        0 & 0
                 \end{pmatrix}  ,0,\cdots
            \right)  \\
    \cdot e^{ 2\pi i
             \left(
                  c_1 \frac{x_{-1}^{1,1}-u_1x_{-1}^{2,1}}{x_{-1}^{2,1}}
                  +\cdots +c_{n-1}\frac{x_{-1}^{n-1,n-1}-u_{n-1}x_{-1}^{n,n-1}}{x_{-1}^{n,n-1}}
                  +c_n \frac{x_{-1}^{n,n}-v_n x_{-2}^{1,n}}{x_{-2}^{1,n}}
             \right)}   dx
 \end{multline*}
which is exactly the right side of (\ref{6.14}).

Therefore the functional (\ref{6.13}) gives a Whittaker functional for the first copy of $\widehat{SL_n}(F)$.
Similarly we will also give the Whittaker functional for the second copy of $\widehat{SL_n}(F)$ as follows:
For $d = (d_1, \cdots, d_n ) \in \Bbb R^n$, define a linear functional $\Psi_d': \cal S(M_n(\Bbb R)[t^{-1}]t^{-1})_\lambda \to \Bbb C$ by
\begin{align} \label{6.15}
\Psi_d'(f)=& \int f
\left(
  \begin{pmatrix}
      x_{-1}^{1,1} & \cdots &x_{-1}^{1,n-1}&x_{-1}^{1,n}\\
      x_{-1}^{2,1} & \cdots & x_{-1}^{2,n-1} & x_{-1}^{2,n} \\
      & \ddots   & \vdots & \vdots\\ 0 && x_{-1}^{n,n-1} & x_{-1}^{n,n}
  \end{pmatrix} ,
  \begin{pmatrix}
      0&\cdots & 0 & x_{-2}^{1,n}\\
      0 & \cdots & 0 & 0\\
      \vdots & \ddots & \vdots & \vdots\\
      0 & \cdots & 0 & 0
  \end{pmatrix},
  0,\cdots
\right)   \nonumber \\
& \cdot e^{2\pi i \left(d_1 \frac{x_{-1}^{2,2}}{x_{-1}^{2,1}} +\cdots +d_{n-1}\frac{x_{-1}^{n,n}}{x_{-1}^{n,n-1}} +
                   d_n \frac{x_{-1}^{1,1}}{x_{-2}^{1,n}}
\right)} dx
\end{align}

Let $\hat\chi_d$ be a character of the unipotent subgroup in $B^T$ defined as follows:
 For any $v\in B^T$ such that $v(0)$ is lower triangular unipotent, that is
  \[
     v=   \begin{pmatrix} 1   & & \\  v_1&  \ddots  & \\ \vdots &   \ddots& \ddots & \\ \ast & \cdots& v_{n-1}&1\end{pmatrix} +
           \begin{pmatrix}   &&&z_n \\ &&&\\ & &&\\&&& \end{pmatrix}t+\cdots,
  \]
 then define
 \[
 \hat \chi_d (v):= e^{-2\pi i(d_1v_1+\cdots+d_{n-1}v_{n-1} +d_n z_{n})},
  \]
 where $v_1,\cdots,v_{n-1}$ are subdiagonal elements of $v\in B^T$, and $z_n$ is the $(1,n)$-entry of $v'(0)$.
Similar to the computation before, we have
 \begin{equation} \label{6.16}
(\pi_2^\ast(v)\Psi_d',  f)=\hat\chi_d(v)(\Psi_d',f).
\end{equation}

\


\begin{thebibliography}{CERP}

\bibitem[ALZ]{ALZ} D. Adamovic, R. Lü, K. Zhao, Whittaker modules for the affine Lie algebra $A_1^{(1)}$,
   Advances in Mathematics,
Volume 289, 2016, Pages 438-479.

\bibitem[Fr]{Fr}  I. B. Frenkel, Two constructions of affine Lie algebra representations and boson-fermion correspondence in quantum field theory, Jour. of Func. Anal.
Vol 44, Issue 3, Pages 259-327, 1981.


 \bibitem[FF]{FF}  A.J. Feingold and I. B. Frenkel, Classical affine algebras. Adv. in Math., 56(2):117 - 172, 1985.


 \bibitem[FZ]{FZ} I.B. Frenkel and A.M.Zeitlin,  On The Continuous Series for $\widehat{sl ( 2 , {\Bbb R})} $.  Communications in Mathematical Physics,
 Vol. 326, Issue 1, pp 145–165, 2014.



\bibitem[GJ]{GJ} R.Godement and H.Jacquet, Zeta Functions of Simple Algebras,
  Lecture Notes in Mathematics,  Vol. 260, Springer-Verlag, Berlin-New York, 1972.

\bibitem[GR]{GR} H.Garland and M.S.Raghunathan, A Bruhat Decomposition for the Loop Space of a Compact Group:
 a New Approach to Results of Bott, Proc. Natl. Acad. Sci.  USA 72 (1975), 4716-4717.



 \bibitem[Ka]{Ka}  V. Kac, Infinite dimensional Lie algebras, 3rd Edition, Cambridge University Press (1990)



 \bibitem[KP]{KP}  V. Kac and D.H. Peterson, Spin and wedge representations of infinite-dimensional Lie algebras and groups, Proc Natl Acad Sci USA. 1981 Jun; 78(6): 3308–3312.



     \bibitem[M]{M}  J.Milnor, On Polylogarithms, Hurwitz Zeta Functions, and the Kubert Identities,    Enseign. Math. (2) 29 (1983), no. 3-4, 281-322.

  \bibitem[Z]{Z}
       Y. Zhu, {\it Theta functions and Weil representations of loop symplectic groups}, Duke Math. J. {\bf 143} (2008), no. 1, 17-39.




\end{thebibliography}
\end{document}